\documentclass[12pt,a4paper]{article}
\RequirePackage{amsmath}%
\RequirePackage{amsthm}%
\usepackage{amsfonts}%
\usepackage{amssymb}%
\usepackage{graphicx}%
\usepackage[autolinebreaks]{mcode}
\usepackage{url}
\usepackage{upquote,textcomp}

\renewcommand{\i}{\mathrm{i}}

\newcommand{\CC}{{\mathbb C}}

\begin{document}
\title{PlgCirMap: A MATLAB toolbox for computing conformal mappings from polygonal multiply connected domains onto circular domains}
\author{Mohamed M.S. Nasser  \\[10pt] 
Department of Mathematics, Statistics and Physics, Qatar University,\\
P.O.Box 2713, Doha, Qatar.\\
{\tt mms.nasser@qu.edu.qa} }
\date{}
\maketitle
\begin{center}
\begin{quotation}
{\noindent {\bf Abstracts.\;\;}%
This paper presents a MATLAB toolbox for computing the conformal mapping from a given polygonal multiply connected domain onto a circular multiply connected domain and its inverse. The toolbox can be used for multiply connected domains with high connectivity and complex geometry. It can be employed also for simply connected domains. 
}%
\end{quotation}
\end{center}
\begin{center}
\begin{quotation}
\end{quotation}
\end{center}
\begin{center}
\begin{quotation}
\end{quotation}
\end{center}

\section{Motivation and significance}
\label{sc:mot}

Conformal mappings are used to transform two-dimensional domains with complex geometry (physical domains) onto domains with simpler one (canonical domains).
Numerous canonical domains have been considered in the literature for conformal mappings of multiply connected domains in the extended complex plane $\CC\cup\{\infty\}$~\cite{koe,Nas-jmaa11,Nas-jmaa13}. 
Perhaps the most important canonical domain for simply and multiply connected domains is the circular domain; i.e., a domain all of whose boundaries are circles. This is due to the existence of analytic formulas for several problems in circular multiply connected domains (see the recent monograph~\cite{Cro-b} and the references cited therein).
Furthermore, circular domains are ideal for using Fourier series and FFT~\cite{Del-Sch04}.

Important examples of complex geometry domains are the polygonal domains, whose boundaries consist of straight line segments.
For simply connected domains, the Schwartz-Christoffel (SC, for short) formula provides us with an explicit form of the conformal mapping from the unit disk onto a given polygonal domain~\cite{Dri02}. 
The SC formula for simply connected domains was discovered independently by Christoffel in 1867 and Schwarz in 1869 (see~\cite[p.~4]{Dri02}).
The generalization of this formula to doubly connected domains was due to Komatu~\cite{Kom} in 1945 (see~\cite[pp.~478-486]{Hen86}). However, the extension of the SC formula to multiply connected domains was established only recently. 
Indeed, DeLillo, Elcrat and Pfaltzgraff~\cite{Del-Sch04} and DeLillo~\cite{Del-Sch06} derived SC formulas for conformal mappings from circular domains onto unbounded and bounded polygonal domains, respectively, using the reflection principle. Crowdy~\cite{Cro-Sch05,Cro-Sch07} presented SC formulas for computing such mappings using Schottky-Klein prime
functions. 

Driscoll~\cite{Dri-ACM96,Dri-ACM05,Dri-sc} created a MATLAB package called SC Toolbox for computing the conformal mapping from the unit disk onto a given polygonal simply connected domain. The toolbox is a generalization of the Fortran package SCPACK developed by Trefethen~\cite{Tre89}. The SC Toolbox has been widely used by many researchers. However, no such toolbox is available so far for polygonal multiply connected domains. The development of such a MATLAB toolbox is the subject of this paper. The proposed toolbox can be used for computing the conformal mapping $w=f(z)$ from a given polygonal multiply connected domain $G$ onto a circular domain $D$ and its inverse $z=f^{-1}(w)$. 

\section{The conformal mapping}

Assume $G$ is a given bounded or unbounded polygonal multiply connected domain bordered by $m$ polygons $\Gamma_j$, $j=1,2,\ldots,m$, such that no corner of these polygons is a cusp. For $m=1$, the domain $G$ is simply connected. If $G$ is bounded, then we assume that $\Gamma_m$ is the external boundary component and encloses all the other boundary components $\Gamma_j$, $j=1,2,\ldots,m-1$. The total boundary $\Gamma=\partial G=\cup_{j=1}^{m}\Gamma_j$ is oriented such that $G$ is on the left of $\Gamma$. 
Then there exists a conformal mapping $w=f(z)$ from the domain $G$ onto a circular multiply connected domain $D$ bordered by $m$ circles $C_j$, $j=1,2,\ldots,m$ (see~\cite[p.~pp. 488-496]{Hen86} and~\cite[pp. 118-127]{Wen}). The mapping $f$ extends to the boundary of $G$ with $f(\Gamma_j)=C_j$, $j=1,2,\ldots,m$. We assume that the circular domain $D$ is bounded when $G$ is bounded and $D$ is unbounded when $G$ is unbounded. The total boundary $C=\partial D=\cup_{j=1}^{m}C_j$ has the same orientation as $\Gamma$. 

When $G$ is bounded, the conformal mapping $f$ is uniquely determined by assuming that the external boundary $C_m=f(\Gamma_m)$ of $D$ is the unit circle and
\begin{equation}\label{eq:cod-b}
f(\alpha)=0, \quad f'(\alpha)>0,
\end{equation} 
where $\alpha$ is a given point in the domain $G$. The condition~(\ref{eq:cod-b}) can be replaced with the condition
\begin{equation}\label{eq:cod-b2}
f(\alpha)=0, \quad f(\beta)=1,
\end{equation} 
 where $\beta$ is a given point on the external boundary $\Gamma_m$. 
On the other hand, if $G$ is unbounded, then the mapping function $f$ is uniquely determined by assuming that 
\begin{equation}\label{eq:cod-u}
f(z)=z+O\left(\frac{1}{z}\right)
\end{equation} 
near infinity. Alternatively, $f$ is uniquely determined by assuming that $C_m=f(\Gamma_m)$ is the unit circle and
\begin{equation}\label{eq:cod-u2}
f(\infty)=\infty, \quad f(\beta)=1,
\end{equation} 
 where $\beta$ is a given point on the boundary $\Gamma_m$.

The SC formulas derived in~\cite{Cro-Sch05,Cro-Sch07,Del-Sch06,Del-Sch04} can be used to compute the inverse mapping $z=f^{-1}(w)$ from the circular domain $D$ onto the polygonal domain $G$. However, using these SC formulas requires solving a system of non-linear equations to determine the preimages of the polygons' vertices as well as the centers and the radius of the circles. Solving such a nonlinear system of equations is still a challenging problem. 
On the other hand, conformal mappings from multiply connected domains onto circular
domains can be computed using Koebe's iterative method~\cite{koe10} (see~\cite[\S~17.7]{Hen86}). As a special case, Koebe's method can be used to compute the conformal mapping $w=f(z)$ from a given polygonal multiply connected domain $G$ onto a circular multiply connected domain $D$. A fast implementation of Koebe's iterative method using the boundary integral equation with the generalized Neumann kernel is given in~\cite{Nas-cmft15} (see also~\cite{Nas-etna15}). The method presented in~\cite{Nas-cmft15} can be used also to compute the inverse map $z=f^{-1}(w)$ from the circular domain $D$ onto the polygonal domain $G$. Further, this method can be applied even for domains with high connectivity and complex geometry, see~\cite{Nas-etna15,Nas-cmft15}. More recently, another efficient numerical method for computing the conformal mapping from circular domains onto polygonal domains based on rational approximations has been presented by Gopal and Trefethen~\cite{GT-cm} and Trefethen~\cite{T-cm}.

In the presented PlgCirMap toolbox, the above described conformal mapping $w=f(z)$ from the polygonal domain $G$ onto the circular domain $D$ as well as its inverse $z=f^{-1}(w)$ from $D$ onto $G$ will be computed using the boundary integral method presented in~\cite{Nas-cmft15}.

\section{Software Framework}
\label{sc:soft}

\subsection{Software Architecture}
\label{sc:soft-arc}

\begin{figure}[ht]%
\centerline{
\scalebox{0.8}{\includegraphics[trim=0.0cm 0.0cm 0.0cm 0.0cm,clip]{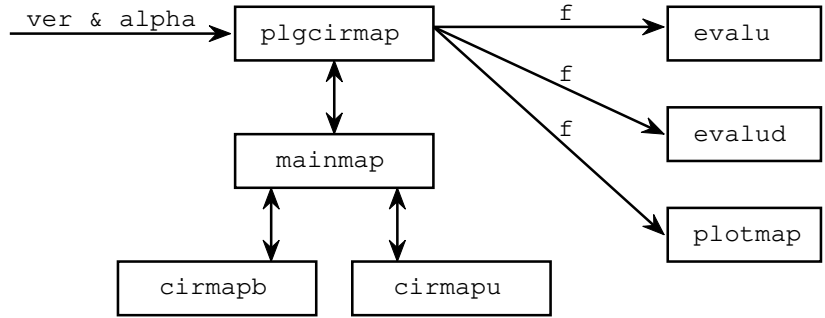}}
}
\caption{PlgCirMap's software architecture.}%
\label{fig:diag}%
\end{figure}

The PlgCirMap is a MATLAB toolbox that consists of different MATLAB functions. The main functions in this toolbox are \verb|plgcirmap|, \verb|evalu|, \verb|evalud|, and \verb|plotmap|. The function \verb|plgcirmap| itself depends on three main functions,  \verb|mainmap|, \verb|cirmapb|, and \verb|cirmapu| (see Figure~\ref{fig:diag}). 
The inputs for the function \verb|plgcirmap| are a cell array \verb|ver| containing the vertices of the polygons and a point \verb|alpha| in the polygonal domain $G$. The default values of the parameters for the numerical calculations are set in the function \verb|plgcirmap|. The numerical implementation of Koebe's iterative method is given in the function \verb|mainmap|. Koebe's iterative method requires computing conformal mappings from bounded simply connected domains onto the unit disk and computing conformal mappings from unbounded simply connected domains onto the exterior of the unit disk (see~\cite{Nas-cmft15} for more details). Such conformal mappings will be computed using the functions \verb|cirmapb| and \verb|cirmapu|, respectively. The output of the \verb|plgcirmap| will be a MATLAB object \verb|f| containing the required information about the conformal mapping $f$ and its inverse $f^{-1}$. From the object \verb|f|, we can compute the values of the conformal mapping and its inverse using the function \verb|evalu|. 
The values of the first derivatives of $f$ and $f^{-1}$ can be computed using the function \verb|evalud|.
Finally, the function \verb|plotmap| uses the object \verb|f| to visualize the conformal mapping $f$ and its inverse $f^{-1}$.

\subsection{Software Functionalities}
\label{sc:soft-fun}

The PlgCirMap toolbox is used to compute and visualize the conformal mapping $w=f(z)$ from a given polygonal multiply connected domain $G$ onto a circular multiply connected domain $D$ and its inverse $z=f^{-1}(w)$. All conditions~(\ref{eq:cod-b})--(\ref{eq:cod-u2}) can be implemented in the toolbox. PlgCirMap can be used also for polygonal simply connected domains. To use the PlgCirMap toolbox, the boundary of the polygonal domain is assumed to have no cusps or slits.

\section{Implementation and Empirical Results}
\label{sc:imp}

\subsection{Parameters' default values}

In the PlgCirMap toolbox, the mapping function $w=f(z)$ and its inverse $z=f^{-1}(w)$ are computed using the boundary integral method presented in~\cite{Nas-cmft15}. The method is based on a fast numerical implementation of Koebe's iterative method using the boundary integral equation with the generalized Neumann kernel~\cite{Nas-etna15,Nas-cmft15}. Assume that each polygon $\Gamma_j$ has $\ell_j\ge3$ vertices. We discretize each segment of the polygon $\Gamma_j$ with $n$ graded mesh points so that $\Gamma_j$ is discretized by $\ell_j n$ graded mesh points (see~\cite{kre,Nas-cmft17} for details on how the graded mesh points are chosen). 
Then, applying the Nystr\"om method with the trapezoidal rule reduces the integral equation to a linear system which is solved iteratively by the MATLAB function \verb|gmres|. 
The matrix-vector product in the GMRES method is computed using the MATLAB function \verb|zfmm2dpart| from the MATLAB toolbox FMMLIB2D~\cite{Gre-Gim12}. 
The computational cost of the method is $O(m^2\ell n+m\ell n\log n)$ where $\ell=\max_{1\le j\le m}\ell_j$ (see~\cite{Nas-etna15,Nas-cmft15} for details). 

In the PlgCirMap toolbox, the default values of the parameters for the numerical calculations are set in the function \verb|plgcirmap| as follows.
\begin{enumerate}
	\item  The default value of $n$, the number of discretization points in each side of the polygons, is set to \verb|n=2^9|. In fact, accurate results can be obtained even for the values of $n$ as small value as \verb|n=2^5|. Increasing the value of $n$ leads to increased accuracy of the obtained results. However, choosing very large values of $n$ should be avoided since it could cause a problem with the convergence of the FMM functions in the toolbox FMMLIB2D~\cite{Gre-Gim12}.

  \item For the FMM function \verb|zfmm2dpart|, we set the default value \verb|iprec=4|, which means the accuracy of the FMM is $0.5\times10^{-12}$. The accuracy of the obtained results can be improved by choosing \verb|iprec=5| (the accuracy of the FMM will be $0.5\times10^{-15}$). However, there may be a problem with the convergence of the FMM if we choose \verb|iprec=5|, especially when $n$ is too large.
		
	\item In the MATLAB function \verb|gmres|, the default tolerances of the GMRES method is set to \verb|gmrestol=0.5e-12| and the default maximum number of iterations allowed is set to \verb|gmresmaxit=100|. The GMRES is used without restart. 

	\item For Koebe's iterative method, the default tolerance and the default maximum number of iterations allowed are set to \verb|koebetol=1e-12| and \verb|koebemaxit=100|, respectively.  		
\end{enumerate}

\section*{Remarks.}
\begin{enumerate}
	\item The default value for maximum number of iterations allowed for both the GMRES method and Koebe's iterative method is $100$. However, for several numerical experiments with the PlgCirMap toolbox, both methods converges with less than $100$ iterations.
	\item If we choose very large values of $n$ and/or \verb|iprec=5|, the FMM functions in the toolbox FMMLIB2D might cause the computer to crash. Unfortunately, no warning message will be displayed and sometimes one need to restart MATLAB. In such case, we need to reduce the value of $n$.
\end{enumerate}

\subsection{The main functions of the toolbox}

\subsubsection{The function {\tt plgcirmap}}
\label{sc:plgcirmap}

To call this function, we need first to define the vertices \verb|ver| of the polygons and a point \verb|alpha| in the domain $G$. 
Here, \verb|ver| is a cell array where \verb|ver{j}| are the vertices of the polygon $\Gamma_j$ for $j=1,2,\ldots,m$. 
The vertices must be ordered such that the domain $G$ is always on the left of the boundary $\Gamma$. 
If $G$ is bounded, the \verb|ver{m}| are the vertices of the external polygon and \verb|alpha| is the point in the domain $G$ that will be mapped to $0$ in the domain $D$. 
When $G$ is unbounded, we define \verb|alpha=inf| and it will be mapped to \verb|inf| in the domain $D$. 

For computing the conformal mapping with the normalizations~(\ref{eq:cod-b}) or~(\ref{eq:cod-u}), the function \verb|plgcirmap| is called as follows:
\[
{\tt f = plgcirmap(ver,alpha)}.
\]
To compute the conformal mapping with the normalizations~(\ref{eq:cod-b2}) or~(\ref{eq:cod-u2}), we call the function \verb|plgcirmap| as:
\[
{\tt f = plgcirmap(ver,alpha,ver\{end\}(k))}
\]
which means that the vertex $k$ on the polygon $\Gamma_m$ will be mapped to $1$.
The object \verb|f|, which is a MATLAB struct with several fields, contains the data of the conformal mapping $w=f(z)$ from $G$ onto $D$ as well as its inverse $z=f^{-1}(w)$ from $D$ onto $G$. For example,
\begin{enumerate}
	\item \verb|f.cent| is a vector of length $m$ containing the centers of the circles $C_j$, $j=1,2,\ldots,m$.
	\item \verb|f.rad| is a vector of length $m$ containing the radius of the circles $C_j$, $j=1,2,\ldots,m$.
	\item \verb|f.imgver| is a cell array where \verb|f.imgver{j}| are the images of the vertices of the polygon $\Gamma_j$ on the circle $C_j$ for $j=1,2,\ldots,m$. Note that these computed values are known in the literature as the \emph{preimages} of the vertices of the polygons.
	\item \verb|f.et| is a vector that contains the discretization of the parametrization of the boundary $\Gamma=\partial G$.
	\item \verb|f.zet| is a vector that contains the discretization of the parametrization of the boundary $C=\partial D$.
\end{enumerate}

\subsubsection{The function {\tt evalu}}
\label{sc:evalu}

Once the function \verb|plgcirmap| is executed and the MATLAB struct \verb|f| is computed, we use the function \verb|evalu| to compute the values of the mapping function $f$ and its inverse $f^{-1}$. 
For computing the direct mapping $f$ at a vector of points $\verb|z|$ in $G$, we call the function $\verb|evalu(f,z,'d')|$.
Similarly, the values of the inverse mapping $f^{-1}$ at a vector of points $\verb|w|$ in $D$ can be computed through $\verb|evalu(f,w,'v')|$.

\subsubsection{The function {\tt evalud}}
\label{sc:evalud}

The function \verb|evalud| is used to compute the values of the first derivatives of the mapping function $f$ and its inverse $f^{-1}$. 
The values of $f'(z)$ at a vector of points $\verb|z|$ in $G$ can be computed via $\verb|evalud(f,z,'d')|$ and the values of $(f^{-1})'(w)$ at a vector of points $\verb|w|$  in $D$ can be computed by $\verb|evalud(f,w,'v')|$.

\subsubsection{The function {\tt plotmap}}
\label{sc:plot}

The function \verb|plotmap| is used to visualize the conformal mapping $f$ and its inverse $f^{-1}$.
To plot rectangular grids in the polygonal domain $G$ and their images in the circular domain $D$, we call $\verb|plotmap(f,'d','rec',n1,n2)|$ where $\verb|n1|$ is the number of horizontal lines and $\verb|n2|$ is the number of vertical lines. We can also plot rectangular grids in $D$ and their images in $G$ through $\verb|plotmap(f,'v','rec',n1,n2)|$.
Similarly, polar grids can be plotted via $\verb|plotmap(f,'d','plr',n1,n2)|$ and $\verb|plotmap(f,'v','plr',n1,n2)|$ where $\verb|n1|$ is the number of circles and $\verb|n2|$ is the number of rays. To plot the domains $G$ and $D$ without grid-points, we call $\verb|plotmap(f)|$.

\subsection{Comparison with Schwarz-Christoffel Toolbox~\cite{Dri-sc}}

\begin{figure}%
\centerline{
\scalebox{0.4}{\includegraphics[trim=1.0cm 0.0cm 1.0cm 0.0cm,clip]{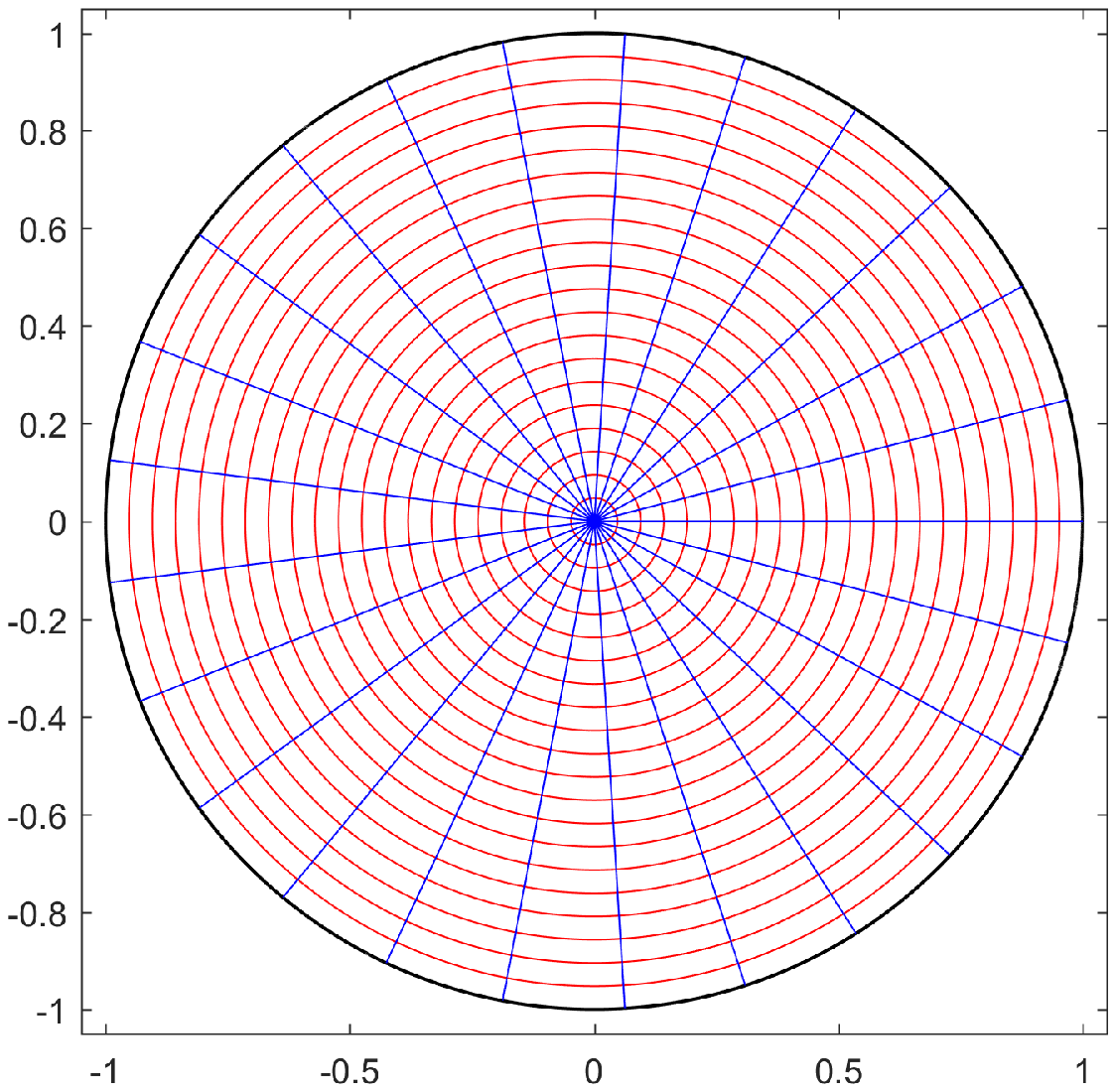}}
\hfill
\scalebox{0.4}{\includegraphics[trim=1.0cm 0.0cm 1.0cm 0.0cm,clip]{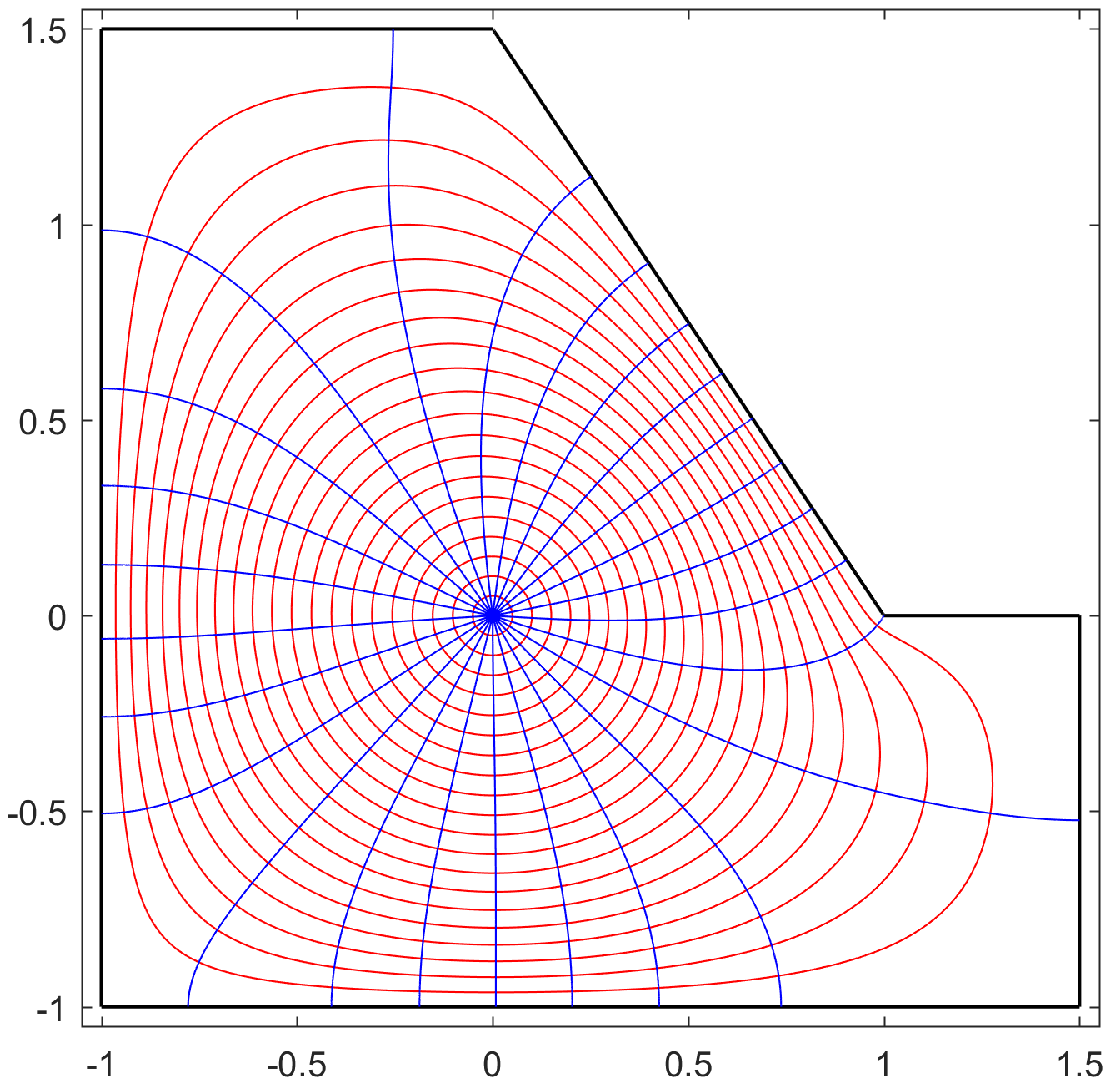}}
\hfill
\scalebox{0.4}{\includegraphics[trim=1.0cm 0.0cm 1.0cm 0.0cm,clip]{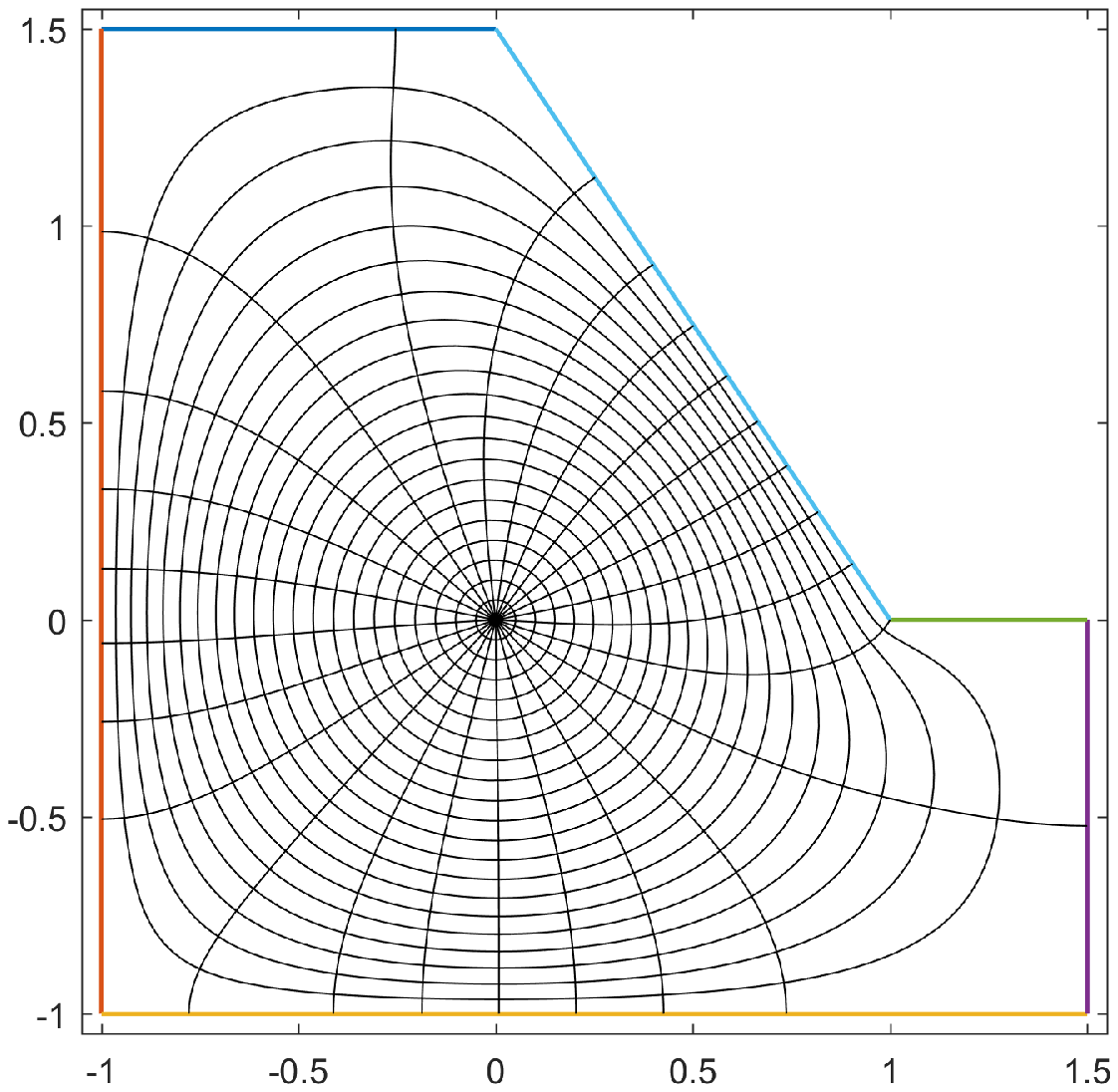}}
}
\caption{A comparison between the PlgCirMap Toolbox (center) and the SC Toolbox (right).}%
\label{fig:sim}%
\end{figure}

The PlgCirMap toolbox can be used for computing the conformal mapping from a given polygonal simply connected domain $G$ onto the unit disk (for bounded $G$) or the exterior of the unit disk (for unbounded $G$). For such case, the accuracy of the PlgCirMap toolbox can be compared against the well known SC Toolbox~\cite{Dri-sc}. Consider the simply connected domain $G$ interior to the polygon with the vertices $1.5\i$, $-1+1.5\i$, $-1-\i$, $1.5-\i$, $1.5$, and $1$ (see Figure~\ref{fig:sim}). The SC Toolbox can be used to map the unit disk onto the domain $G$ such that $1$ on the unit circle is mapped to the last vertex of the polygonal which is also $1$. So, we shall assume here that $f$ is the conformal mapping from the polygonal domain $G$ onto the unit disk $D$ normalized by~(\ref{eq:cod-b2}) with $\alpha=0$ and $\beta=1$, i.e., $f(0)=0$ and $f(1)=1$.

To use the PlgCirMap toolbox in computing such conformal mapping $f$, we first set the vertices of the polygon and the point $\alpha$ in $G$ as follows:
\begin{verbatim}
>> ver{1}=[1.5i ; -1+1.5i ; -1-1i ; 1.5-1i ; 1.5 ; 1];
>> alpha = 0;
\end{verbatim}
Then, we use the MATLAB function \verb|plgcirmap| to compute the object \verb|f|,
\begin{verbatim}
>> f=plgcirmap(ver,alpha,ver{1}(6));
\end{verbatim}
To plot orthogonal polar grids in the unit disk $D$ (see Figure~\ref{fig:sim} (left)) and their images under the inverse mapping $f^{-1}$ in the polygonal domains $G$, we call
\begin{verbatim}
>> plotmap(f,'v','plr',20,25); 
\end{verbatim}
The resulting figure is shown in Figure~\ref{fig:sim} (center).

Next, we use the SC Toolbox to compute the inverse mapping $f^{-1}$ from the unit disk $D$ onto the polygonal domain $G$. First, we set the accuracy of the SC toolbox to $10^{-14}$ by calling
\begin{verbatim}
>> options = scmapopt('Tolerance',1e-14);
\end{verbatim}
Then, we use the SC toolbox to compute a MATLAB object \verb|fisc| by calling
\begin{verbatim}
>> p=polygon(ver{1});
>> fisc = diskmap(p,options); 
>> fisc = center(fisc,alpha);
\end{verbatim}
The plot of the image of the orthogonal polar grids in the disk $D$ under the inverse mapping $f^{-1}$ computed by the SC toolbox is shown in Figure~\ref{fig:sim} (right). This figure is generated by calling
\begin{verbatim}
>> plot(fisc,20,25);
\end{verbatim}

For comparison, we compute the preimages of the vertices of the polygon using the PlgCirMap toolbox by calling
\begin{verbatim}
>> prevertpcm = f.imgver{:};
\end{verbatim}
and using the SC toolbox by calling
\begin{verbatim}
>> prevertsc  = get(fisc,'prevert');
\end{verbatim}
Then, we compute the maximum norm between the computed values as
\begin{verbatim}
>> E_1 = norm(prevertsc-prevertpcm,inf)
\end{verbatim}
The obtained maximum norm is $E_1=1.0579\times10^{-12}$.

Next, we choose a set of points 
\begin{verbatim}
>> zz  = 0.6.*exp(i.*linspace(0,2*pi,1000));
\end{verbatim}
in the polygonal domain $G$. To compute the values of the conformal mapping $f$ at these points \verb|zz| using the SC toolbox, we call the function \verb|evalinv(fisc,zz)|. For the PlgCirMap toolbox, the values of the mapping $f$ at the points \verb|zz| are computed by calling \verb|evalu(f,zz,'d')|. Then, we compute the maximum norm between the computed values as
\begin{verbatim}
>> E_2 = norm(evalinv(fisc,zz)-evalu(f,zz,'d'),inf)
\end{verbatim}
The obtained maximum norm is $E_2=5.3952\times10^{-13}$.

We also compare the two toolboxes by choosing a set of points 
\begin{verbatim}
>> ww = 0.9.*exp(i.*linspace(0,2*pi,1000));
\end{verbatim}
in the unit disk $D$. Then, we compute the maximum norm between the values of the inverse map $f^{-1}$ computed at the points \verb|ww| by the two toolboxes as
\begin{verbatim}
>> E_3 = norm(fisc(ww)-evalu(f,ww,'v'),inf)
\end{verbatim}
The obtained maximum norm is $E_3=1.0291\times10^{-12}$.

Finally, we check the accuracy of each toolbox separately. We use both toolboxes to compute approximate values for $f^{-1}(f({\tt zz}))$ and $f(f^{-1}({\tt ww}))$. Then, we compute the maximum error norm in the computed values for the SC toolbox by
\begin{verbatim}
>> ES_1 = norm(fisc(evalinv(fisc,zz))-zz,inf)
>> ES_2 = norm(evalinv(fisc,fisc(ww))-ww,inf)
\end{verbatim}
The obtained values are $ES_1=9.6538\times10^{-15}$ and $ES_2=1.4457\times10^{-14}$.
For the PlgCirMap toolbox, we compute the maximum error norm in the computed values by
\begin{verbatim}
>> E_4 = norm(evalu(f,evalu(f,zz,'d'),'v')-zz,inf)
>> E_5 = norm(evalu(f,evalu(f,ww,'v'),'d')-ww,inf)
\end{verbatim}
and the obtained values are $E_4=7.3621\times10^{-14}$ and $E_5=1.1931\times10^{-13}$.

As we see from this example, there is a very good agreement between the results obtained by the well developed SC toolbox and the presented PlgCirMap toolbox (for the default value of $n$, $n=2^9$). Further, we will have a good agreement between the results obtained by the two toolboxes even for small values of $n$. Indeed, if we change the value of $n$ in the MATLAB function \verb|plgcirmap| to $n=2^5$, then the above computed norms will be as follows: $E_1=9.6102\times10^{-8}$, $E_2=1.0283\times10^{-8}$, $E_3=6.3408\times10^{-8}$, $E_4=4.7186\times10^{-9}$, and $E_5=8.2561\times10^{-9}$.

\section{Illustrative Examples}
\label{sc:exam}

In this section, we use the PlgCirMap toolbox to compute the conformal mappings for two domains. More examples for both bounded and unbounded domains are available in \url{https://github.com/mmsnasser/PlgCirMap}. 

\subsection{A bounded multiply connected domain}
\begin{figure}%
\centerline{
\hfill
\scalebox{0.4}{\includegraphics[trim=1.0cm 0.0cm 1.0cm 0.0cm,clip]{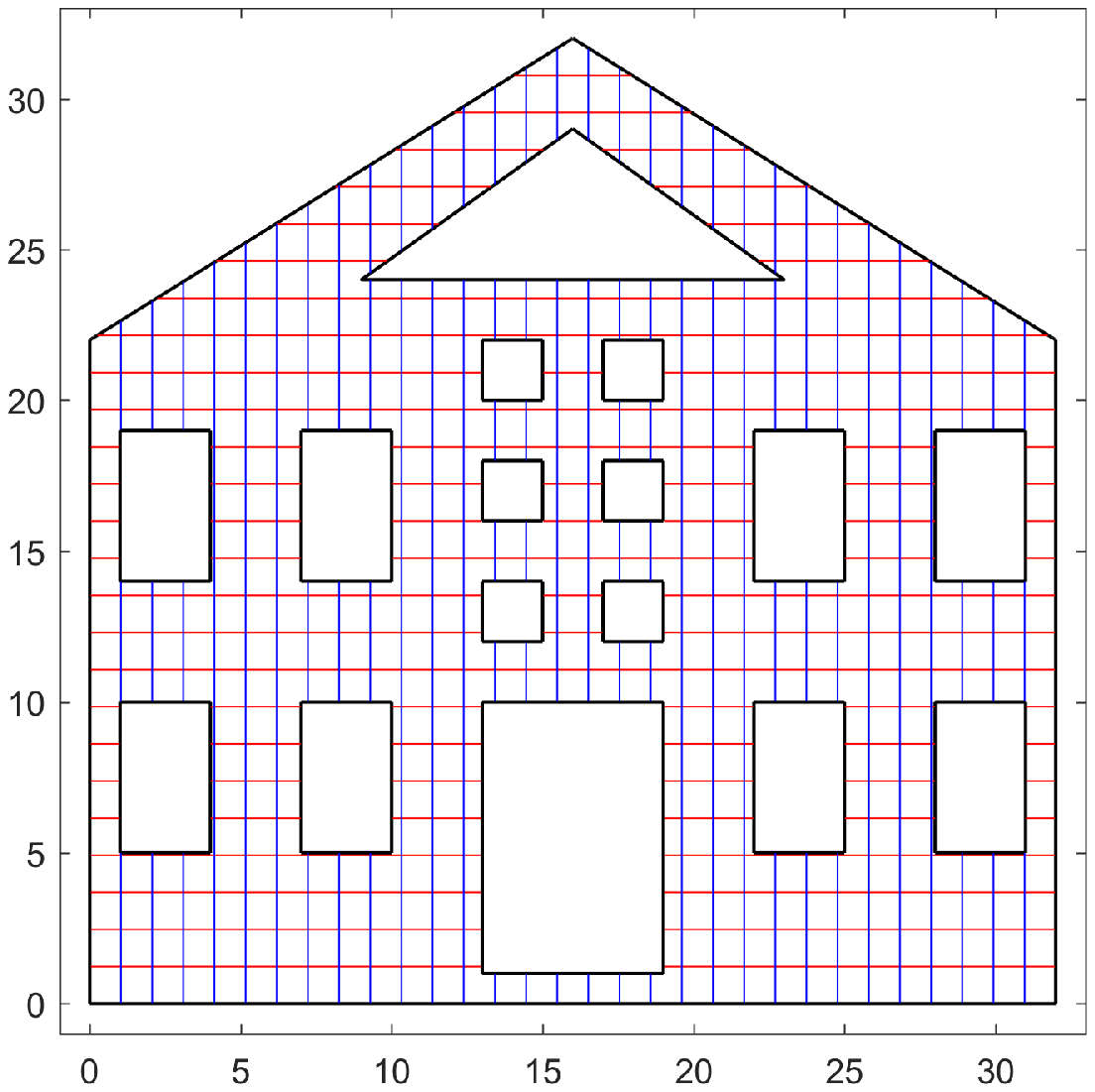}}
\hfill
\scalebox{0.4}{\includegraphics[trim=1.0cm 0.0cm 1.0cm 0.0cm,clip]{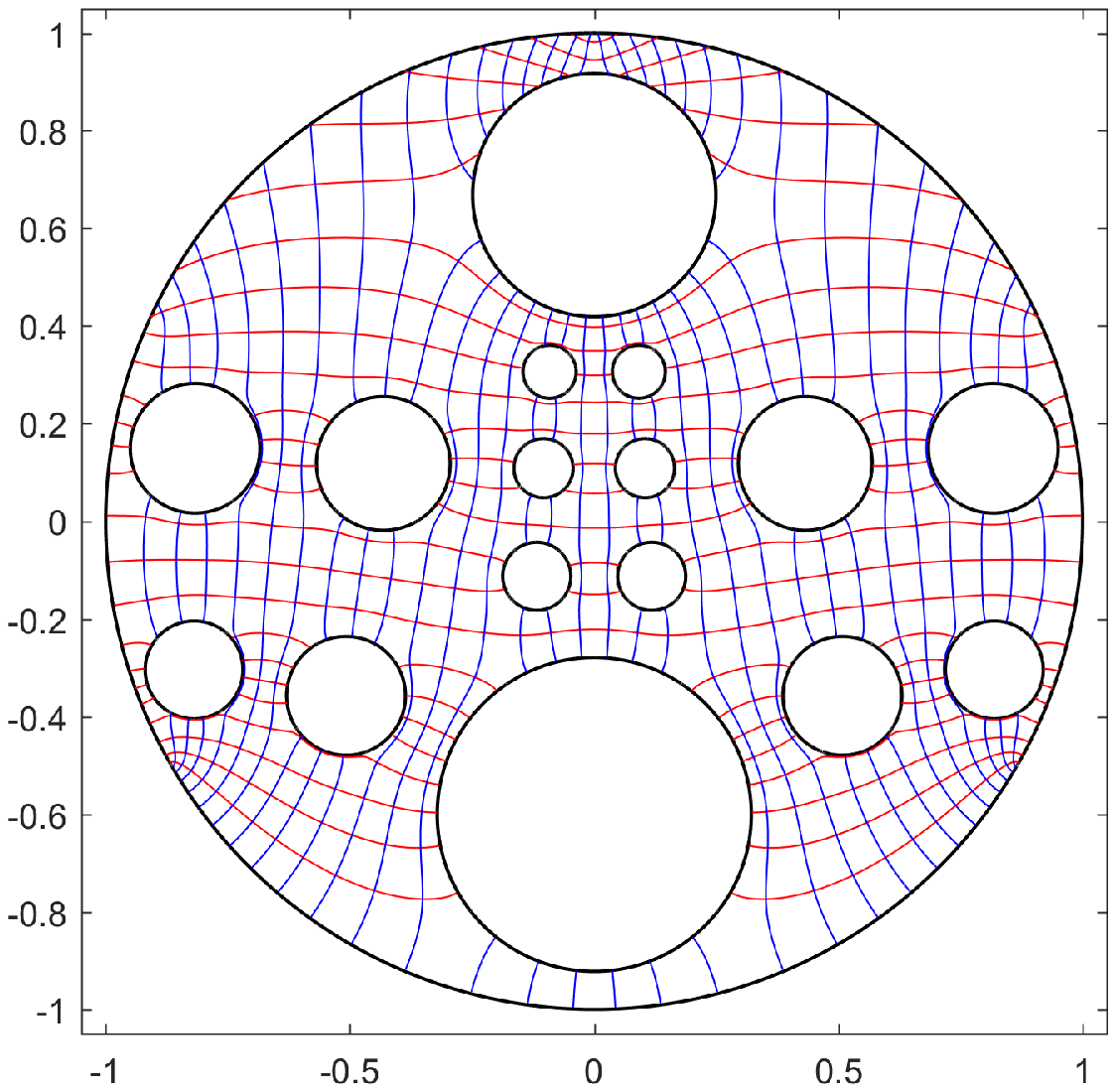}}
\hfill
}
\centerline{
\hfill
\scalebox{0.4}{\includegraphics[trim=1.0cm 0.0cm 1.0cm 0.0cm,clip]{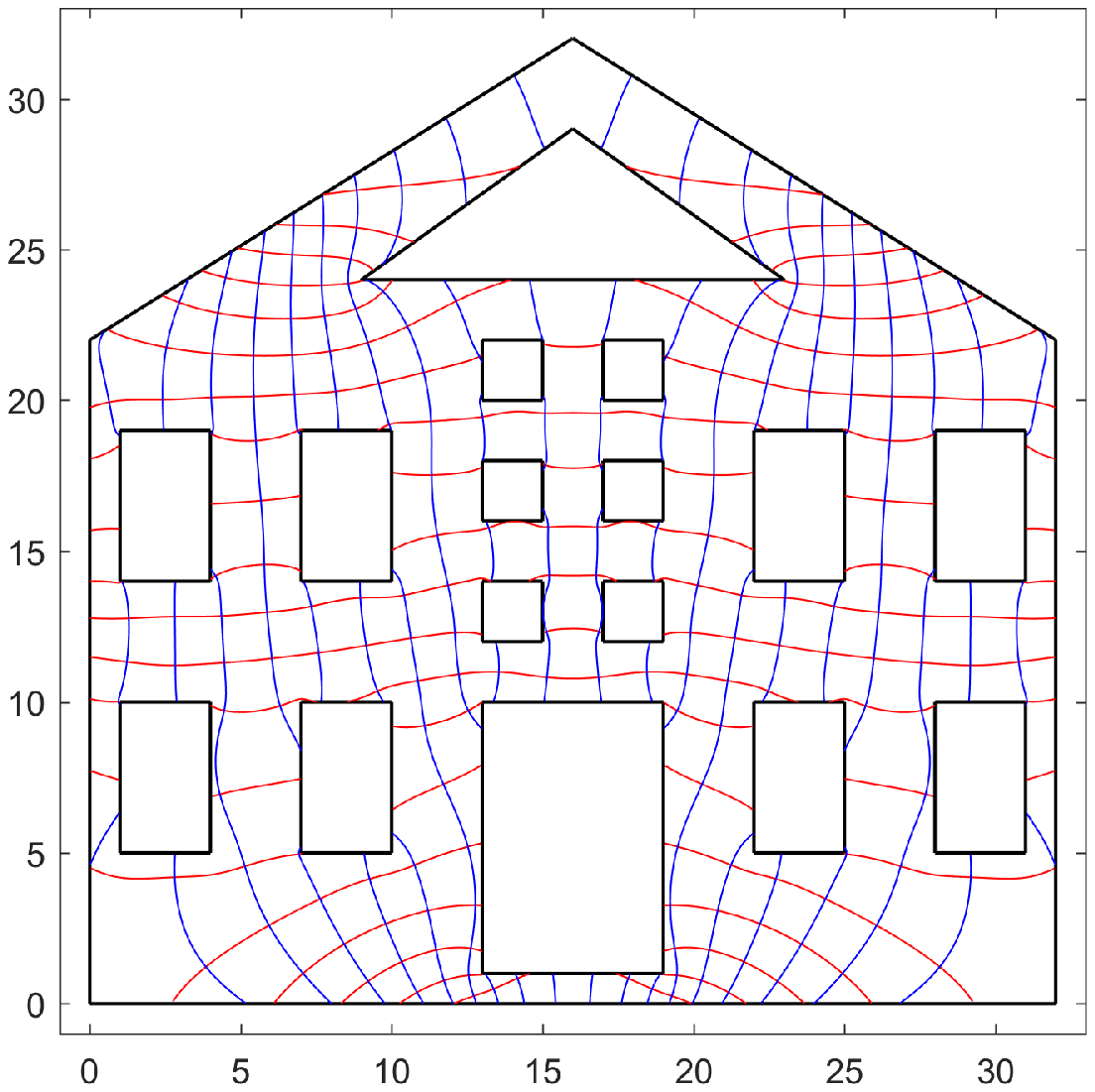}}
\hfill
\scalebox{0.4}{\includegraphics[trim=1.0cm 0.0cm 1.0cm 0.0cm,clip]{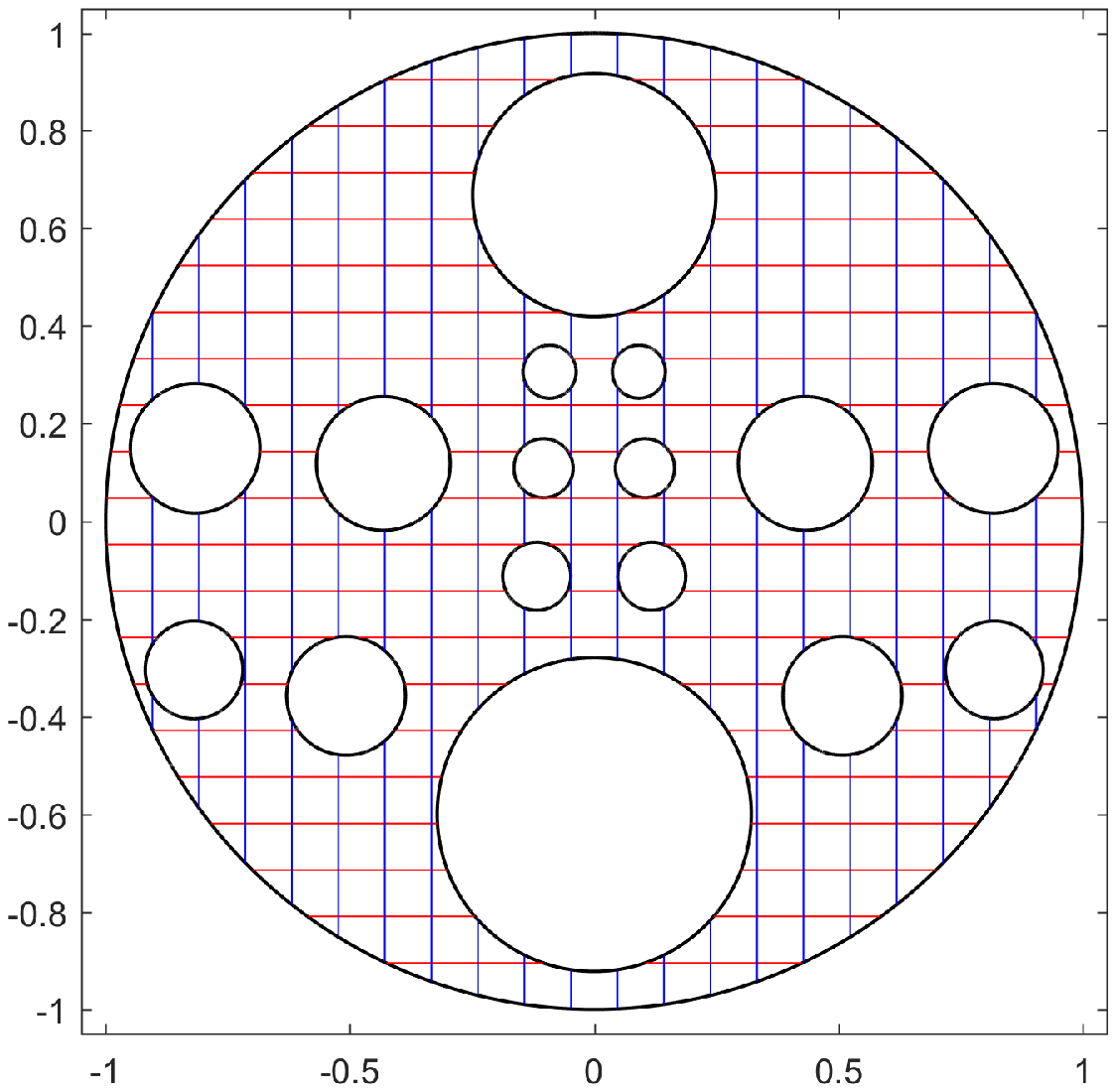}}
\hfill
}
\centerline{
\hfill
\scalebox{0.4}{\includegraphics[trim=1.0cm 0.0cm 1.0cm 0.0cm,clip]{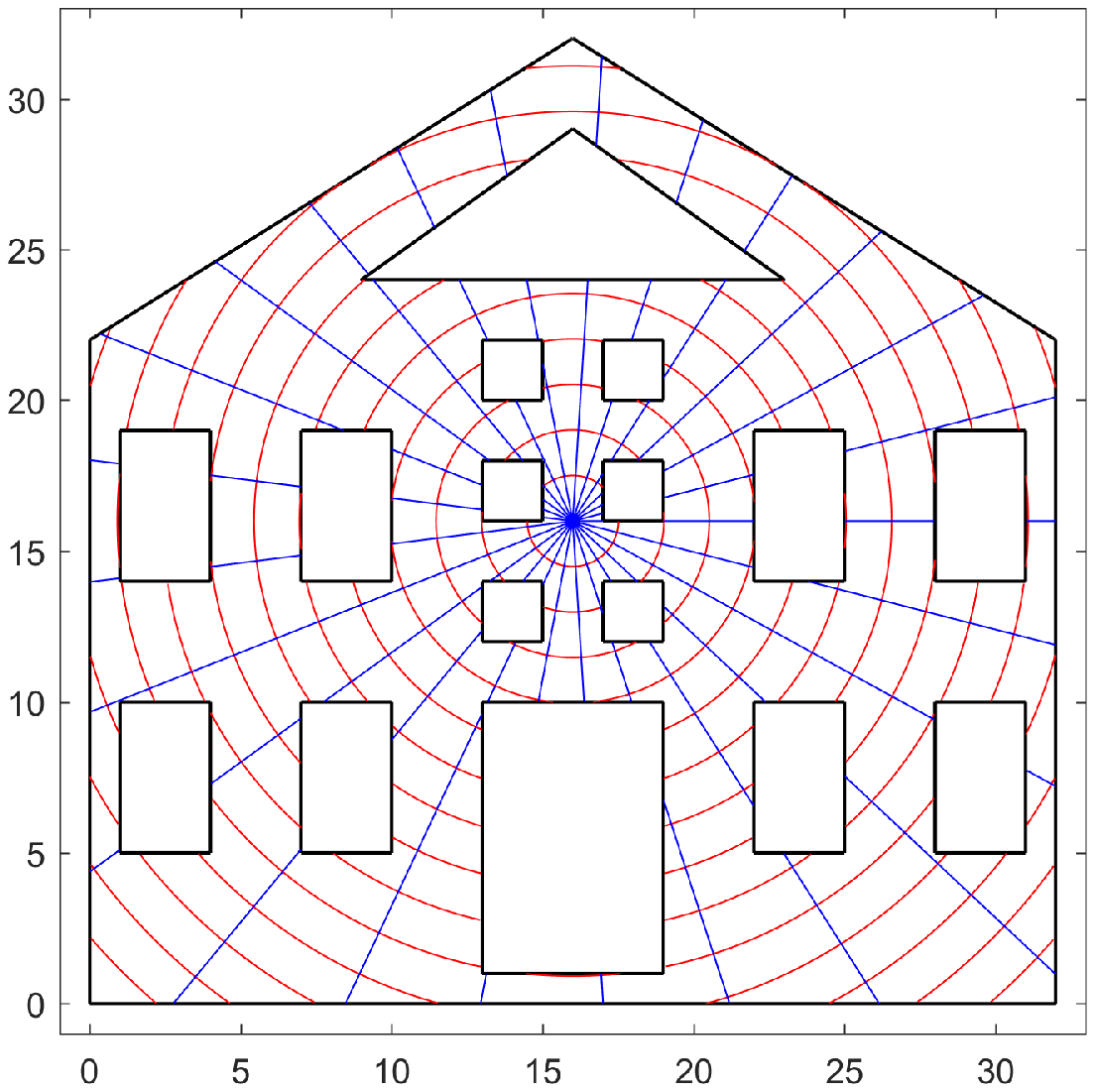}}
\hfill
\scalebox{0.4}{\includegraphics[trim=1.0cm 0.0cm 1.0cm 0.0cm,clip]{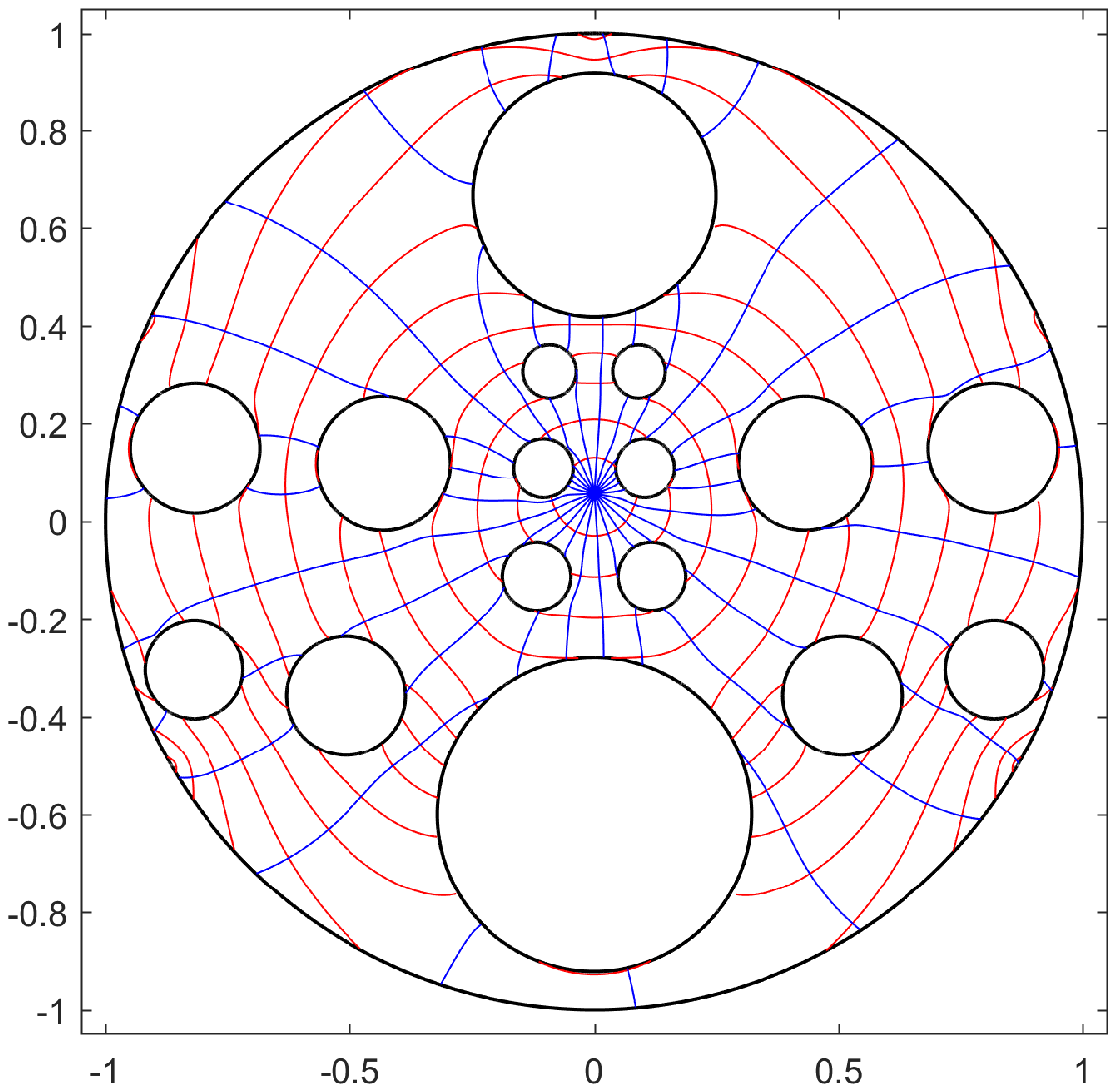}}
\hfill
}
\centerline{
\hfill
\scalebox{0.4}{\includegraphics[trim=1.0cm 0.0cm 1.0cm 0.0cm,clip]{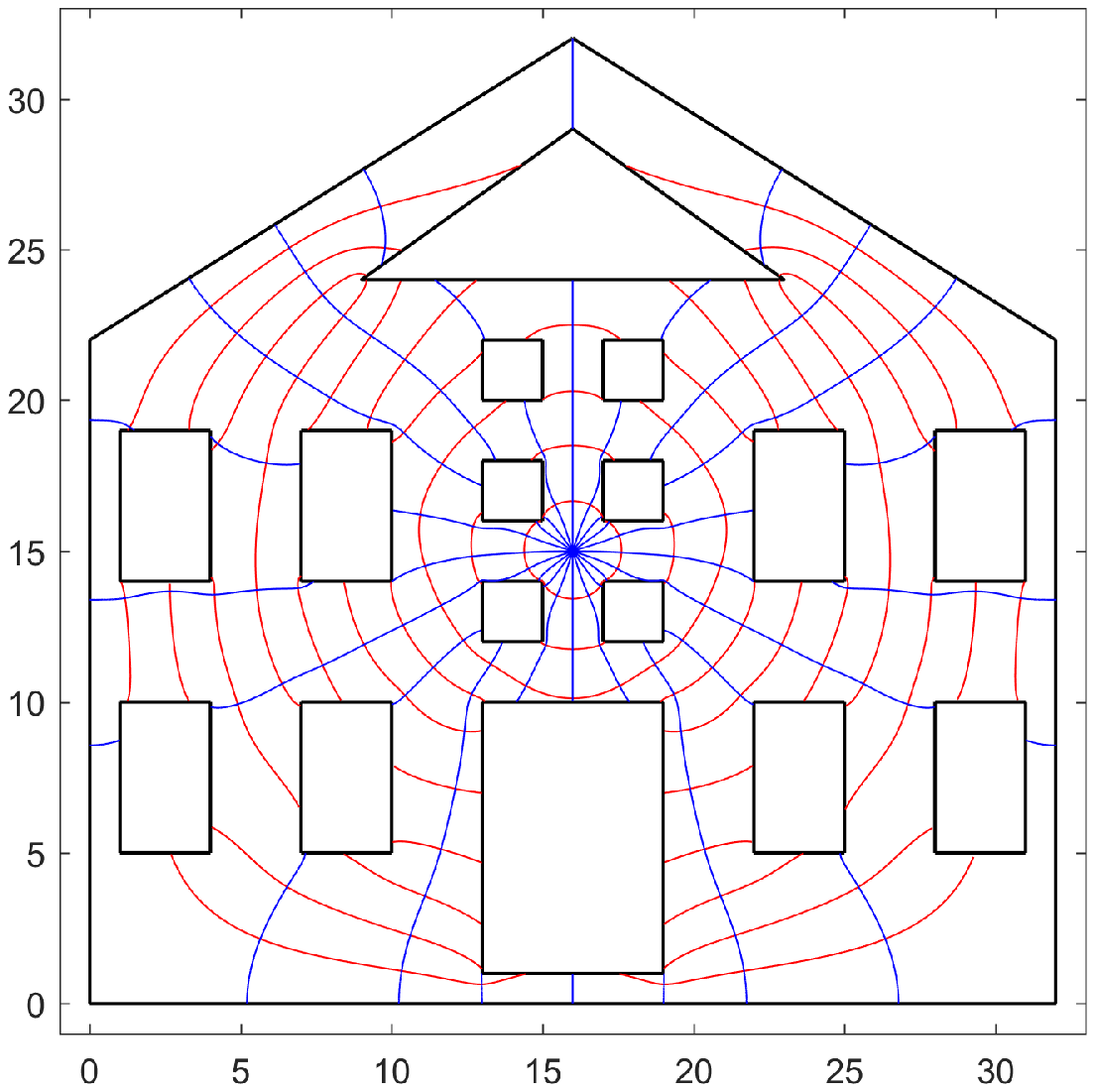}}
\hfill
\scalebox{0.4}{\includegraphics[trim=1.0cm 0.0cm 1.0cm 0.0cm,clip]{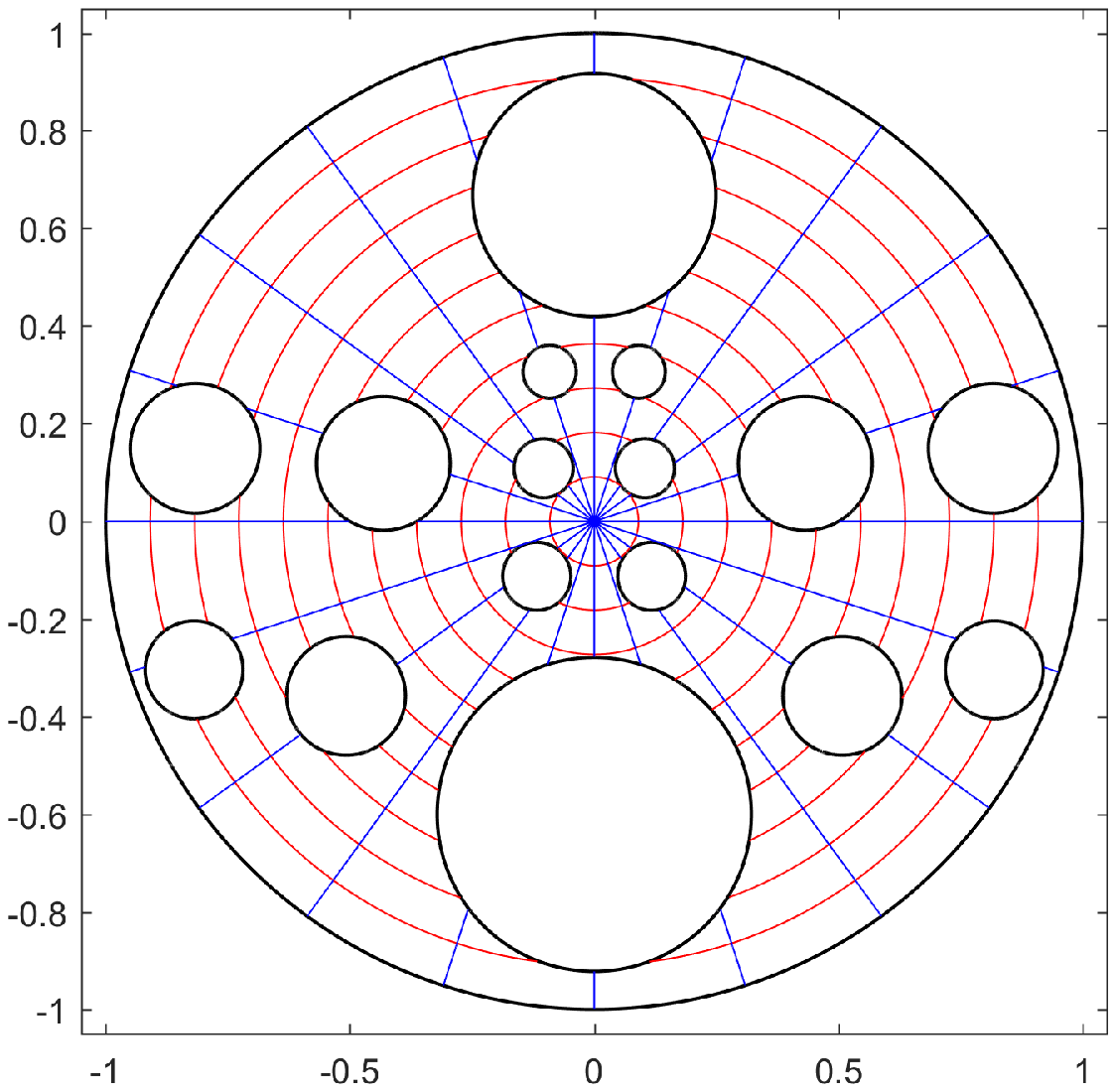}}
\hfill
}
\caption{An example of a bounded multiply connected domain.}%
\label{fig:mult}%
\end{figure}

For the first example, we consider a bounded polygonal multiply connected domain of connectivity $17$. In the MATLAB code given below, we first define the vertices of the polygons and we choose a point $\alpha$ in the domain $G$. Then, we call the function $\verb|plgcirmap|$ to compute the conformal mapping with the normalization~(\ref{eq:cod-b}). 
The function $\verb|plotmap|$ is then called to visualize the conformal mapping $f$ and its inverse $f^{-1}$ as in Figure~\ref{fig:mult}.

\begin{lstlisting}
ver{1}  = [31+10i ; 31+5i  ; 28+5i  ; 28+10i ];
ver{2}  = [25+10i ; 25+5i  ; 22+5i  ; 22+10i ];
ver{3}  = [19+10i ; 19+1i  ; 13+1i  ; 13+10i ];
ver{4}  = [10+10i ; 10+5i  ;  7+5i  ;  7+10i ];
ver{5}  = [ 4+10i ;  4+5i  ;  1+5i  ;  1+10i ];
ver{6}  = [31+19i ; 31+14i ; 28+14i ; 28+19i ];
ver{7}  = [25+19i ; 25+14i ; 22+14i ; 22+19i ];
ver{8}  = [19+14i ; 19+12i ; 17+12i ; 17+14i ];
ver{9}  = [15+14i ; 15+12i ; 13+12i ; 13+14i ];
ver{10} = [19+18i ; 19+16i ; 17+16i ; 17+18i ];
ver{11} = [15+18i ; 15+16i ; 13+16i ; 13+18i ];
ver{12} = [19+22i ; 19+20i ; 17+20i ; 17+22i ];
ver{13} = [15+22i ; 15+20i ; 13+20i ; 13+22i ];
ver{14} = [10+19i ; 10+14i ;  7+14i ;  7+19i ];
ver{15} = [ 4+19i ;  4+14i ;  1+14i ;  1+19i ];
ver{16} = [16+29i ; 23+24i ;  9+24i ];
ver{17} = [16+32i ;  0+22i ;  0+0i  ; 32+0i  ; 32+22i ];
alpha = 16+15i;
f = plgcirmap(ver,alpha);
plotmap(f,'d','rec',25,30); 
plotmap(f,'v','rec',20,20); 
plotmap(f,'d','plr',15,25); 
plotmap(f,'v','plr',10,20); 
\end{lstlisting}

\subsection{An unbounded multiply connected domain}
\begin{figure}%
\centerline{
\scalebox{0.5}{\includegraphics[trim=0.0cm 0.0cm 0.0cm 0.0cm,clip]{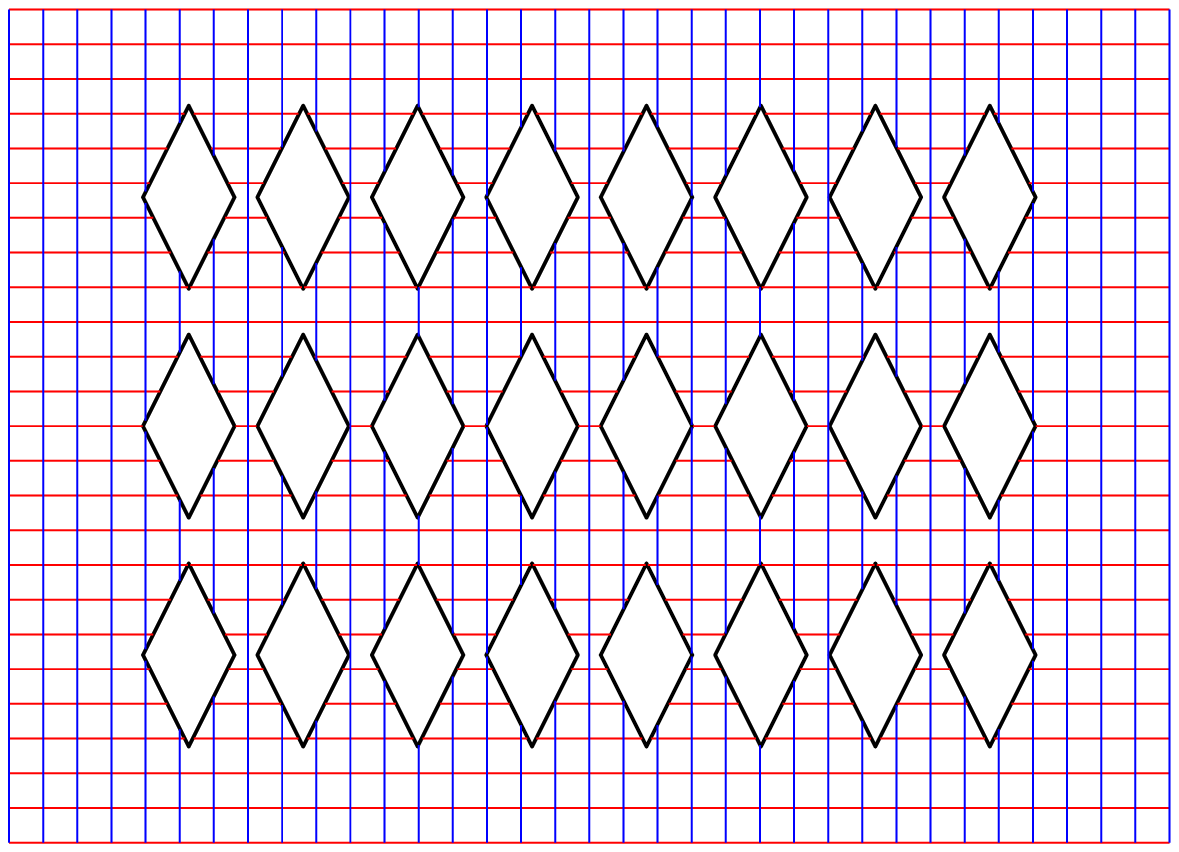}}
\hfill
\scalebox{0.5}{\includegraphics[trim=0.0cm 0.0cm 0.0cm 0.0cm,clip]{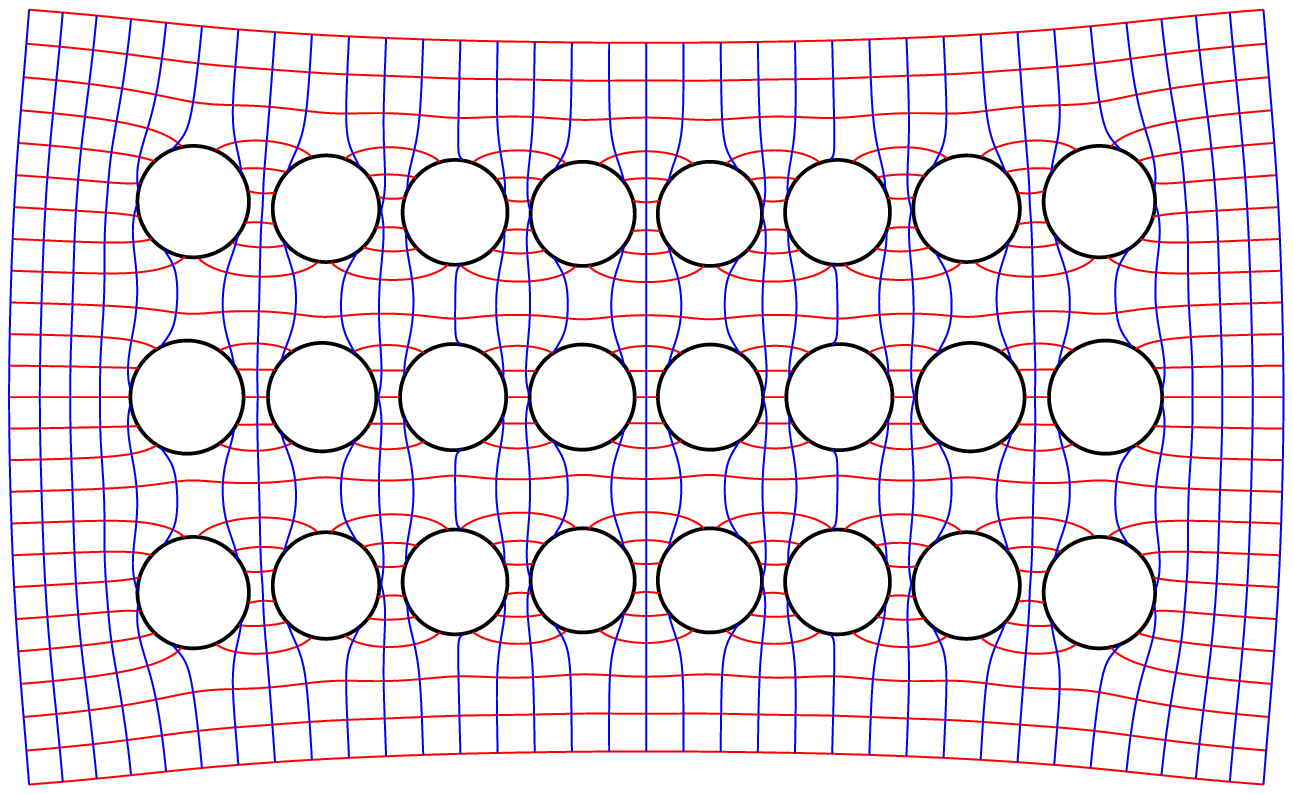}}
}
\caption{An example of an unbounded multiply connected domain.}%
\label{fig:multu}%
\end{figure}

In the second example, we consider an unbounded polygonal multiply connected domain of connectivity $24$. 
The MATLAB code for this example is given below where the normalization~(\ref{eq:cod-u}) is used. Figure~\ref{fig:multu} shows the obtained figures.

\begin{lstlisting}
ver{1}  = [-1.75-0.9i ;-1.95-0.5i ;-1.75-0.1i ;-1.55-0.5i ];
ver{2}  = [-1.25-0.9i ;-1.45-0.5i ;-1.25-0.1i ;-1.05-0.5i ];
ver{3}  = [-0.75-0.9i ;-0.95-0.5i ;-0.75-0.1i ;-0.55-0.5i ];
ver{4}  = [-0.25-0.9i ;-0.45-0.5i ;-0.25-0.1i ;-0.05-0.5i ];
ver{5}  = [ 0.25-0.9i ; 0.05-0.5i ; 0.25-0.1i ; 0.45-0.5i ];
ver{6}  = [ 0.75-0.9i ; 0.55-0.5i ; 0.75-0.1i ; 0.95-0.5i ];
ver{7}  = [ 1.25-0.9i ; 1.05-0.5i ; 1.25-0.1i ; 1.45-0.5i ];
ver{8}  = [ 1.75-0.9i ; 1.55-0.5i ; 1.75-0.1i ; 1.95-0.5i ];
ver{9}  = [-1.75+0.1i ;-1.95+0.5i ;-1.75+0.9i ;-1.55+0.5i ];
ver{10} = [-1.25+0.1i ;-1.45+0.5i ;-1.25+0.9i ;-1.05+0.5i ];
ver{11} = [-0.75+0.1i ;-0.95+0.5i ;-0.75+0.9i ;-0.55+0.5i ];
ver{12} = [-0.25+0.1i ;-0.45+0.5i ;-0.25+0.9i ;-0.05+0.5i ];
ver{13} = [ 0.25+0.1i ; 0.05+0.5i ; 0.25+0.9i ; 0.45+0.5i ];
ver{14} = [ 0.75+0.1i ; 0.55+0.5i ; 0.75+0.9i ; 0.95+0.5i ];
ver{15} = [ 1.25+0.1i ; 1.05+0.5i ; 1.25+0.9i ; 1.45+0.5i ];
ver{16} = [ 1.75+0.1i ; 1.55+0.5i ; 1.75+0.9i ; 1.95+0.5i ];
ver{17} = [-1.75+1.1i ;-1.95+1.5i ;-1.75+1.9i ;-1.55+1.5i ];
ver{18} = [-1.25+1.1i ;-1.45+1.5i ;-1.25+1.9i ;-1.05+1.5i ];
ver{19} = [-0.75+1.1i ;-0.95+1.5i ;-0.75+1.9i ;-0.55+1.5i ];
ver{20} = [-0.25+1.1i ;-0.45+1.5i ;-0.25+1.9i ;-0.05+1.5i ];
ver{21} = [ 0.25+1.1i ; 0.05+1.5i ; 0.25+1.9i ; 0.45+1.5i ];
ver{22} = [ 0.75+1.1i ; 0.55+1.5i ; 0.75+1.9i ; 0.95+1.5i ];
ver{23} = [ 1.25+1.1i ; 1.05+1.5i ; 1.25+1.9i ; 1.45+1.5i ];
ver{24} = [ 1.75+1.1i ; 1.55+1.5i ; 1.75+1.9i ; 1.95+1.5i ];
alpha = inf;
f=plgcirmap(ver,alpha);
plotmap(f,'d','rec',25,25);
\end{lstlisting}

\section{Impact}
\label{sec:impact}

Conformal mappings are a powerful tool to solve several problems in the fields of science and engineering involving the Laplace equation due to its invariant under conformal mappings. 
With the help of conformal mappings, solving the Laplace equation in domains with complex geometry (physical domains) is reduced to solving this equation in domains with simpler geometry (canonical domains). The simple geometry of the circular domain makes it an important canonical domain from both physical and computational points of view. For example, D. Crowdy with several collaborators have recently presented analytic formulas for several problems of fluid mechanics in circular multiply connected domains (see e.g.,~\cite{Cro-b} and the references cited therein). With the help of the presented toolbox, such analytic formulas can be used also for polygonal domains (we refer to~\cite{Kal} for an example of such applications of the toolbox).

\section{Final remarks and suggestions for future improvements}
\label{sec:rem}

\begin{enumerate}
	\item The numerical method used in the PlgCirMap toolbox to compute the conformal mapping is based on the boundary integral method presented in~\cite{Nas-cmft15}. Thus, the accuracy of the toolbox depend on the accuracy of the numerical solution of the used boundary integral equation. In the current version of the toolbox, the integral equation is solved using the Nystr\"om method with the trapezoidal rule based on graded mesh points (see~\cite{kre} for details). Improving the accuracy of the numerical solution of the integral equation will improve the accuracy of the presented toolbox. This could be considered in the future.
	\item The method presented in~\cite{Nas-cmft15} can be used for general multiply connected domains with smooth or piecewise smooth boundaries. As a result, the PlgCirMap toolbox can be generalized to multiply connected domains other than polygonal domains. In particular, this toolbox can be generalized easily to compute the conformal mapping from multiply connected domains with circular-arc boundaries onto circular multiply connected domains.
	\item The accuracy of the numerical methods for computing the conformal mapping from a circular domain $D$ onto an elongated domain $G$ is seriously affected by what is known as the \emph{crowding phenomenon} (see~\cite[\S2.6]{Dri02}). Although the method used in this toolbox is based on computing the conformal mapping from the elongated domain $G$ onto the circular domain $D$, it still affected by crowding. One possible way to improve the accuracy of the presented toolbox for elongated domains in future is to use the domain decomposition method or to consider canonical domains other than the circular domain (see~\cite{Dri02}).

\end{enumerate}

\section*{Acknowledgements}

This toolbox has been presented in the workshop: ``The complex analysis toolbox: new techniques and perspectives'' which held in ``Isaac Newton Institute for Mathematical Sciences'', Cambridge, September 9-13, 2019. The author would like to thank INI for the hospitality during the workshop.

\end{document}